\documentclass[11pt]{amsart}
\usepackage[mathscr]{eucal}
\usepackage{amsfonts}
\usepackage{amsmath}
\usepackage{amsthm}
\usepackage{amssymb}
\usepackage{latexsym}
\theoremstyle{plain}

\newtheorem{Thm}[equation]{Theorem}
\newtheorem{Cor}[equation]{Corollary}
\newtheorem{Prop}[equation]{Proposition}
\newtheorem{Lem}[equation]{Lemma}

\numberwithin{equation}{section}

\newcommand{\bx}{\Box}
\newcommand{\iso}{\stackrel{\sim}{\rightarrow}}
\newcommand{\la}{\langle}
\newcommand{\ra}{\rangle}
\renewcommand{\u}{\underline}

\newcommand{\A}{{\mathscr{A}}}
\newcommand{\BB}{{\mathscr{B}}}
\newcommand{\Bi}{{\mathscr{B}i\mathscr{L}ie}}
\newcommand{\B}{{\mathbf{B}}}
\newcommand{\C}{{\mathbf{C}}} 
\newcommand{\D}{{\mathbf{D}}}
\newcommand{\Det}{{\textrm{Det}}}
\newcommand{\edge}{{\textrm{edge}}}
\newcommand{\Edge}{{\textrm{Edge}}}
\newcommand{\E}{{\mathscr{E}nd}}
\newcommand{\F}{{\mathscr{F}}}
\newcommand{\Hom}{{\mathscr{H}om}}
\newcommand{\I}{{\mathscr{I}}}
\newcommand{\Id}{{\mathtt{Id}}}
\newcommand{\In}{{\textrm{In}}}
\newcommand{\Ind}{{\mathtt{Ind}}}
\newcommand{\K}{{\mathbf{K}}}
\newcommand{\Ker}{{\mathtt{Ker}}}
\newcommand{\LL}{{\mathbf{L}}}
\newcommand{\op}{{\mathtt{op}}}
\newcommand{\Out}{{\textrm{Out}}}
\newcommand{\PP}{{\mathscr{P}}}
\newcommand{\QQ}{{\mathscr{Q}}}
\newcommand{\sgn}{{\textrm{Sgn}}}
\newcommand{\Ver}{{\textrm{Vert}}}

\def\k{{\Bbbk}}
\def\s{{\mathbb{S}}}

\begin{document}
\title{KOSZUL DUALITY FOR DIOPERADS}

\author{Wee Liang Gan}
\address{Department of Mathematics, University of Chicago, Chicago,
IL 60637, U.S.A.}
\email{wlgan@math.uchicago.edu}

\begin{abstract}
We introduce the notion of a dioperad to describe certain operations 
with multiple inputs and multiple outputs. The framework of 
Koszul duality for operads is generalized to dioperads.
We show that the Lie bialgebra dioperad is Koszul.
\end{abstract}

\maketitle

\section*{\bf Introduction}

The current interest in the understanding of various algebraic
structures using operads is partly due to the theory of Koszul duality
for operads, which was developed by Ginzburg-Kapranov in \cite{GiK}; see
eg. \cite{Ka} or \cite{L} for surveys.
However, algebraic structures such as bialgebras and Lie bialgebras,
which involve both multiplication and comultiplication, or bracket 
and cobracket, are defined using PROP's (cf. \cite {Ad}) rather than
operads.
Inspired by the theory of string topology of Chas-Sullivan (cf. 
\cite{ChS}, \cite{Ch}, \cite{Tr}), Victor Ginzburg suggested to the author
that there should be a theory of Koszul duality for PROP's. 
This paper results from the attempt to develop such a theory.
More precisely, we introduce the notion of a dioperad, which can be used
to describe certain operations with multiple inputs and multiple outputs.
We show that one can set up a theory of Koszul duality for dioperads,
and we prove that the dioperad associated to Lie bialgebra
is Koszul.

Let us explain how dioperads arise. 
Suppose $\QQ = \{\QQ(m,n)\}$ is a 
PROP defined by some generators and relations.
We think of an element in $\QQ(m,n)$ as an operation
obtained by compositions of the generators according to
a ``flow chart'' with $n$ inputs and $m$ outputs.
Assuming that the defining relations
between the generators of $\QQ$ are expressed by flow charts 
which are trees, there is a subspace $\PP(m,n)$ of $\QQ(m,n)$ 
consisting of those operations obtained from the flow charts
which are trees.
The collection $\PP = \{\PP(m,n)\}$ is precisely the
dioperad with the same generators and relations as $\QQ$.
Since the defining relations of $\QQ$ are expressed by trees,
no essential information is lost by restricting our
attention to $\PP$.
If $f$ and $g$ are operations
in $\PP$, then by substituition of the $j$-th output of $g$ into
the $i$-th input of $f$, we get another operation 
$f_{i}\circ_{j}g$ which is still in $\PP$. Pictorially, this 
just means that if we join a root of a tree to a leaf of 
another tree, we still get a tree.

The assumption that the defining relations  
between the generators of $\QQ$ are expressed by trees
is satisfied, for example, by the Lie bialgebra PROP, but not
by the bialgebra PROP; see eg. \cite{ES} p.19 and p.70.
Thus, Lie bialgebras can be understood in the framework of
dioperads, but bialgebras cannot.
The above assumption is motivated by the question of what is a
``quadratic'' PROP. It seems an answer 
(suggested by the fact that the differential in Kontsevich's
graph complex \cite{Ko} is induced by edge contractions)
is that the defining relations between the generators 
should be express by graphs with precisely one internal edge; cf.
also \cite{Ge} \S4.8. 
However, a graph with one internal edge is necessarily a tree. 

The paper is organized as follows.
In \S\ref{sec:1}, we give the definition of a dioperad and other
generalities.
In \S\ref{sec:2}, we define the notion of a quadratic dioperad,
its quadratic dual, and introduce our main example of Lie bialgebra
dioperad. In \S\ref{sec:3}, we define the cobar dual of a dioperad.
A quadratic dioperad is Koszul if its cobar dual is quasi-isomorphic
to its quadratic dual.
In \S\ref{sec:4}, we prove a proposition to be used later in
\S\ref{sec:5}. This proposition is a generalization of a result of
Shnider-Van Osdol \cite{SVO}. In \S\ref{sec:5}, we prove that 
Koszulity of a quadratic dioperad is equivalent to exactness of 
certain Koszul complexes. In the case of operads, this is
due to Ginzburg-Kapranov, with a different proof by Shnider-Van Osdol.
The Koszulity of the Lie bialgebra dioperad
follows from this and an adaptation of results of Markl \cite{M2}.

\section{\bf Dioperads}\label{sec:1}

\subsection{}\label{sec:11}
We give in this subsection the definition of dioperad which is 
similar to the definition of operad in \cite{M1}; cf. also \cite{M}.

Let $\mathcal{C}$ be the symmetric monoidal
category of finite dimensional 
differential $\mathbb{Z}$-graded super vector spaces over a field $\k$
of characteristic $0$, and let $\Hom$ be the internal hom functor 
of $\mathcal{C}$.
Let $\s_{n}$ denote the automorphism group of $\{1, \ldots, n\}$.
If $m=m_{1}+\cdots+m_{n}$ is an ordered partition and 
$\sigma \in \s_{n}$, then the \emph{block permutation} $\sigma_{m_{1},
\ldots, m_{n}} \in \s_{m}$ is the permutation that acts on $\{1, \ldots,
m\}$ by permuting $n$ intervals of lengths $m_{1}, \ldots, m_{n}$ in the
same way that $\sigma$ permutes $1, \ldots, n$.
If $\sigma_{1} \in \s_{n_{1}}$, $\sigma_{2} \in \s_{n_{2}}$ and
$i \in \{1, \ldots, n_{1}\}$, then
define $\sigma_{1} \circ_{i} \sigma{_2} \in \s_{n_{1}+n_{2}-1}$
by
\[ \sigma_{1} \circ_{i} \sigma{_2} := (\sigma_{1})_{1, \ldots, 1, n_{2}, 1
\ldots, 1} \circ 
(\Id \times \cdots \times \sigma_{2} \times \cdots \times \Id), \]
where $\sigma_{2}$ is at the $i$-th place. (See eg. \cite{MSS}
Definition 1.2 or \cite{SVO} p.387. )

A \emph{dioperad} $\PP$ in $\mathcal{C}$ consists of data:\\
(i) objects $\PP(m,n)$, equipped with a $(\s_{m},\s_{n})$-bimodule
structure, for each ordered pair of positive integers $m,n$;\\
(ii) morphisms $_{i}\circ_{j}: \PP(m_{1}, n_{1}) \otimes
\PP(m_{2}, n_{2}) \rightarrow \PP(m_{1}+m_{2}-1, n_{1}+n_{2}-1)$ for each
$m_{1}, m_{2}, n_{1}, n_{2} \geq 1$ and $1 \leq i \leq n_{1}$,
$1 \leq j \leq m_{2}$;\\
(iii) a morphism $\eta: \k \rightarrow \PP(1,1)$ such that
\[ _{1}\circ_{i}(\eta \otimes \Id): \k \otimes \PP(m,n)
\iso \PP(m,n) \]
and
\[ _{j}\circ_{1}(\Id \otimes \eta): \PP(m,n) \otimes \k
\iso \PP(m,n) \]
are the canonical isomorphisms for all $m,n \geq 1$ and
$1 \leq i \leq m$, $1 \leq j \leq n$.

These data are required to satisfy the following associativity and
equivariance conditions:\\
(a) for all $m_{1}, n_{1}, m_{2}, n_{2}, m_{3}, n_{3} \geq 1$ and
$1 \leq i \leq n_{1}+n_{2}-1$, $1 \leq j \leq m_{3}$,
$1 \leq k \leq n_{1}$, $1 \leq l \leq m_{2}$, the morphism
\[
_{i}\circ_{j}(_{k}\circ_{l} \otimes \Id):
\PP(m_{1},n_{1}) \otimes \PP(m_{2},n_{2}) \otimes \PP(m_{3},n_{3})
\rightarrow \PP(m_{1}+m_{2}+m_{3}-2 ,n_{1}+n_{2}+n_{3}-2)
\]
is equal to
\[
\left\{ \begin{array}{ll}
(\sigma, 1)
(_{k+n_{3}-1}\circ_{l}) (_{i}\circ_{j} \otimes \Id) (\Id \otimes \tau)
& \textrm{if $i \leq k-1$}\\
_{k}\circ_{j+l-1} (\Id \otimes _{i-k+1}\circ_{j})
& \textrm{if $k \leq i \leq k+n_{2}-1$}\\
(\sigma, 1)
(_{k}\circ_{l}) (_{i-n_{2}+1}\circ_{j} \otimes \Id) (\Id \otimes \tau)
& \textrm{if $k+n_{2} \leq i$}
\end{array} \right.
\]
where \[ \tau: \PP(m_{2},n_{2}) \otimes \PP(m_{3},n_{3}) \iso
\PP(m_{3},n_{3}) \otimes \PP(m_{2},n_{2})\] is the symmetry
isomorphism, and $\sigma \in \s_{m_{1}+m_{2}+m_{3}-2}$ is the block
permutation \[ ((12)(45))_{l-1,j-1,m_{1},m_{3}-j,m_{2}-l} . \] \\
(b) for all $m_{1}, n_{1}, m_{2}, n_{2}, m_{3}, n_{3} \geq 1$ and
$1 \leq i \leq n_{1}$, $1 \leq j \leq m_{2}+m_{3}-1$,
$1 \leq k \leq n_{2}$, $1 \leq l \leq m_{3}$, the morphism
\[
_{i}\circ_{j}(\Id \otimes _{k}\circ_{l}):
\PP(m_{1},n_{1}) \otimes \PP(m_{2},n_{2}) \otimes \PP(m_{3},n_{3}) 
\rightarrow \PP(m_{1}+m_{2}+m_{3}-2 ,n_{1}+n_{2}+n_{3}-2)
\]
is equal to
\[
\left\{ \begin{array}{ll}
(1, \sigma)
(_{k}\circ_{l+m_{1}-1}) (\Id \otimes _{i}\circ_{j}) (\tau \otimes \Id)
& \textrm{if $j \leq l-1$}\\
_{k+i-1}\circ_{l} (_{i}\circ_{j-l+1} \otimes \Id)
& \textrm{if $l \leq j \leq l+m_{2}-1$}\\
(1, \sigma)
(_{k}\circ_{l}) (\Id \otimes _{i}\circ_{j-m_{2}+1}) (\tau \otimes \Id)
& \textrm{if $l+m_{2} \leq j$}\\
\end{array} \right.
\]
where \[ \tau: \PP(m_{1},n_{1}) \otimes \PP(m_{2},n_{2}) \iso
\PP(m_{2},n_{2}) \otimes \PP(m_{1},n_{1})\] is the symmetry
isomorphism, and $\sigma \in \s_{n_{1}+n_{2}+n_{3}-2}$ is the block
permutation \[ ((12)(45))_{i-1,k-1,n_{3},n_{2}-k,n_{1}-i} . \] \\
(c) for all $m_{1}, n_{1}, m_{2}, n_{2} \geq 1$,
$1 \leq i \leq n_{1}$, $1 \leq j \leq m_{2}$, and $\pi_{1} \in
\s_{m_{1}}$, $\sigma_{1} \in \s_{n_{1}}$, $\pi_{2} \in \s_{m_{2}}$,
$\sigma_{2} \in \s_{n_{2}}$,
the morphism
\[ 
_{i}\circ_{j}((\pi_{1}, \sigma_{1}) \otimes (\pi_{2}, \sigma_{2})):
\PP(m_{1},n_{1}) \otimes \PP(m_{2},n_{2}) \rightarrow \PP(m_{1}+m_{2}-1,
n_{1}+n_{2}-1)
\]
is equal to
\[ (\pi_{2}\circ_{\pi_{2}^{-1}(j)}\pi_{1}, 
\sigma_{1}\circ_{i}\sigma_{2})
(_{\sigma_{1}(i)}\circ_{\pi_{2}^{-1}(j)}). \]
{\bf Remarks.} (1) We shall regard an operad $\PP = \{\PP(n)\}$ as a 
dioperad via $\PP(1,n):= \PP(n)$ and $\PP(m,n):= 0$ for $m > 1$.\\
(2) It may be of interests to relate dioperads to
cyclic operads \cite{GeK}
or pseudo-tensor categories \cite{BD},
but the author does not know how to do this.

\subsection{}
Let $V$ be an object of $\mathcal{C}$. Then the \emph{endomorphism
dioperad} $\E_{V}$ is defined by
\[ \E_{V}(m,n) = \Hom(V^{\otimes n}, V^{\otimes m}).\]
If $f\in\E_{V}(m_{1},n_{1}), g\in\E_{V}(m_{2},n_{2})$, then
$f_{i}\circ_{j}g\in\E_{V}(m_{1}+m_{2}-1,n_{1}+n_{2}-1)$ is the morphism
\[ (\Id \otimes \cdots \otimes f \otimes \cdots \otimes \Id) \sigma
(\Id \otimes \cdots \otimes g \otimes \cdots \otimes \Id):
V^{\otimes(n_{1}+n_{2}-1)} \rightarrow V^{\otimes(m_{1}+m_{2}-1)},\]
where $f$ is at the $j$-th place, $g$ is at the $i$-th place,
and $\sigma \in \s_{n_{1}+m_{2}-1}$ is the block permutation
$((12)(45))_{i-1, j-1, 1, m_{2}-j, n_{1}-i}$.

A \emph{morphism} $f: \PP \rightarrow \QQ$ of dioperads in $\mathcal{C}$
is a collection of morphisms $f(m,n): \PP(m,n) \rightarrow \QQ(m,n)$, $m,n
\geq 1$, compatible with the structures of dioperads. If $\PP$ is a
dioperad in $\mathcal{C}$ and $V$ is an object of $\mathcal{C}$ equipped
with a morphism $\PP \rightarrow \E_{V}$, then $V$ is called a
\emph{$\PP$-algebra}.

\subsection{}
If $\PP$ is a dioperad, then its \emph{opposite}
$\PP^{\op}(m,n) := \PP(n,m)$ with the transposed actions is a
dioperad via the composition rule \[ _{i} \stackrel{\op}{\circ}_{j} := 
(_{j}\circ_{i})\tau : \PP^{\op}(m_{1},n_{1}) \otimes
\PP^{\op}(m_{2},n_{2})
\rightarrow \PP^{\op}(m_{1}+m_{2}-1,n_{1}+n_{2}-1),\] where
\[ \tau: \PP(n_{1},m_{1}) \otimes \PP(n_{2},m_{2}) \iso
\PP(n_{2},m_{2}) \otimes \PP(n_{1},m_{1})\] is the symmetry
isomorphism.

If $\PP$ and $\QQ$ are dioperads, then
$(\PP\otimes\QQ)(m,n) := \PP(m,n) \otimes \QQ(m,n)$ with the diagonal 
actions of $\s_{m}$ and $\s_{n}$ is a dioperad via the composition rule
{\setlength\arraycolsep{2pt} \begin{eqnarray*}
_{i}\stackrel{\PP\otimes\QQ}{\circ}_{j} :=
(_{i}\stackrel{\PP}{\circ}_{j} \otimes _{i}\stackrel{\QQ}{\circ}_{j})
(\Id\otimes\tau\otimes\Id) & : & (\PP\otimes\QQ)(m_{1},n_{1}) \otimes
(\PP\otimes\QQ)(m_{2},n_{2}) {} \\ & & {}
\rightarrow (\PP\otimes\QQ)(m_{1}+m_{2}-1,n_{1}+n_{2}-1),
\end{eqnarray*} where
\[ \tau: \QQ(m_{1},n_{1}) \otimes \PP(m_{2},n_{2}) \iso
\PP(m_{2},n_{2}) \otimes \QQ(m_{1},n_{1}) \]
is the symmetry isomorphism.

\subsection{}
The \emph{suspension dioperad} $\Sigma$ is the endomorphism dioperad of
$\k[1]$, and the \emph{desuspension dioperad} $\Sigma^{-1}$ is the   
endomorphism dioperad of $\k[-1]$; see eg. \cite{Sm} \S 3.1. Thus,
$\Sigma(m,n)$ is a one dimensional vector space placed in degree $n-m$
with sign representations of $\s_{m}$ and $\s_{n}$, and $\Sigma^{-1}(m,n)$
is a one dimensional vector space placed in degree $m-n$ with sign
representations of $\s_{m}$ and $\s_{n}$. The suspension of a dioperad
$\PP$ is $\Sigma\PP:= \Sigma\otimes\PP$, and the desuspension is
$\Sigma^{-1}\PP:= \Sigma^{-1}\otimes\PP$.

Observe that if $\PP$ is a dioperad, then
$\{\PP(m,n)[2m-2]\}$ is also a dioperad.
Define the \emph{sheared suspension} dioperad $\Lambda$ by $\Lambda(m,n)
:= \Sigma(m,n)[2-2m]$, and the  \emph{sheared desuspension} dioperad
$\Lambda^{-1}$ by $\Lambda^{-1}(m,n) := \Sigma^{-1}(m,n)[2m-2]$. 
Thus, $\Lambda(m,n)$ is a one dimensional vector space placed in degree
$m+n-2$ with sign representations of $\s_{m}$ and $\s_{n}$, and
$\Lambda^{-1}(m,n)$ is a one dimensional vector space placed in degree
$2-m-n$ with sign representations of $\s_{m}$ and $\s_{n}$.
The sheared suspension of a dioperad $\PP$ is $\Lambda\PP:=
\Lambda\otimes\PP$, and the sheared desuspension is $\Lambda^{-1}\PP:=
\Lambda^{-1}\otimes\PP$.

\section{\bf Quadratic dual}\label{sec:2}
{}From now on, for all dioperads $\PP$, assume $\PP(1,1)=\k$ unless
otherwise stated.

\subsection{}
By a \emph{tree}, we shall always mean a directed tree such that each 
vertex has at least one outgoing edge and at least one incoming edge. 
A \emph{forest} is a disjoint union of trees. 
If each vertex of a forest has valency at least three, then the forest 
is said to be \emph{reduced}. 
Denote by $\Ver(T)$ (resp. $\edge(T)$, $\Edge(T)$)
the set of vertices (resp. internal edges, all edges) of a forest $T$, and
$\Out(v)$ (resp. $\In(v))$ the set of outgoing (resp. incoming) edges
at a vertex $v$. Let $\det(T) := \bigwedge^{|\edge(T)|}
\k^{\edge(T)}$ and $\Det(T) := \bigwedge^{|\Edge(T)|}
\k^{\Edge(T)}$.

A \emph{$(m,n)$-tree} is defined to be a tree with leaves labelled by
$\{1, \dots, n\}$ and roots labelled by $\{1, \dots, m\}$.
For any $(m,n)$-tree $T$, we have the formula
\begin{eqnarray}\label{eqn:2f}
\sum_{v\in\Ver(T)} (|\Out(v)|+|\In(v)|-2) = m+n-2.
\end{eqnarray}
Observe that the maximal number of vertices
in a reduced $(m,n)$-tree is $m+n-2$. 
Note also that a trivalent $(m,n)$-tree has
$m-1$ vertices with two outgoing edges and
$n-1$ vertices with two incoming edges.

If $\u{m}= (m_{1}, \ldots, m_{k})$, where $m_{1}, \ldots, m_{k}$ are
positive integers, then let $|\u{m}|:= m_{1}+\cdots+m_{k}$.
Suppose  $\u{m} = (m_{1}, \ldots, m_{k})$ and $\u{n} = (n_{1}, \ldots,
n_{k})$. A \emph{$(\u{m},\u{n})$-forest} is defined to be 
a disjoint union of trees $T_{1}, \ldots, T_{k}$ such that, for each $i$,
the leaves of $T_{i}$ are labelled by
\[\{n_{1}+\cdots+n_{i-1}+1, \ldots, n_{1}+\cdots+n_{i-1}+n_{i}\}\]
and the roots of $T_{i}$ are labelled by
\[\{m_{1}+\cdots+m_{i-1}+1, \ldots, m_{1}+\cdots+m_{i-1}+m_{i}\}.\]

\subsection{}
Let $E(m,n)$, $m,n \geq1$, be a collection of finite dimensional
$(\s_{m}, \s_{n})$-bimodules with $E(1,1)=0$.
If $T$ is a tree, let
\[ E(T) = \bigotimes_{v \in \Ver(T)} E(\Out(v), \In(v)). \]
The \emph{free dioperad} $\F(E)$ generated by $E$ is defined by
\[ \F(E)(m,n) := \bigoplus_{(m,n)-\textrm{trees $T$}} E(T). \]

Let $(\Lambda E)(m,n):= \Lambda(m,n)\otimes E(m,n)$ and 
$(\Lambda^{-1} E)(m,n):= \Lambda^{-1}(m,n)\otimes E(m,n)$.
Observe that the canonical injections $(\Lambda E)(m,n) \hookrightarrow
(\Lambda\F(E))(m,n)$ and $(\Lambda^{-1} E)(m,n) \hookrightarrow 
(\Lambda^{-1}\F(E))(m,n)$ induce canonical dioperad isomorphisms
\begin{eqnarray}\label{eqn:2a} 
\F(\Lambda E) \iso \Lambda\F(E) \quad\textrm{and}\quad
\F(\Lambda^{-1} E) \iso \Lambda^{-1}\F(E)
\end{eqnarray} respectively.

\subsection{}
Let $\PP$ be a dioperad. An \emph{ideal} $\I$ in $\PP$ is
a collection of $(\s_{m}, \s_{n})$-sub-bimodules $\I(m,n) \subset
\PP(m,n)$ such that $f _{i}\circ_{j} g \in \I(m_{1}+m_{2}-1,
n_{1}+n_{2}-1)$ for all $1 \leq i \leq n_{1}$, $1 \leq j \leq m_{2}$,
whenever $f \in \I(m_{1},n_{1}), g \in \PP(m_{2}, n_{2})$,
or $f \in \PP(m_{1},n_{1}), g \in \I(m_{2}, n_{2})$.
Intersection of ideals is an ideal. If $\I$ is an ideal in $\PP$, then
$(\PP/ \I)(m,n) := \PP(m,n) / \I(m,n)$ is a dioperad.
The augmentation ideal of $\PP$ is defined by
$\bar{\PP}(1,1) := 0$ and $\bar{\PP}(m,n) := \PP(m,n)$ if $m+n \geq 3$.

\subsection{}
Let $E(1,2)$ be a right $\s_{2}$-module, $E(2,1)$ be a left   
$\s_{2}$-module, and $E(m,n) = 0$ for 
$(m,n) \neq (1,2), (2,1)$. Observe that
\[ \F(E)(1,3)= \Ind_{\s_{2}}^{\s_{3}} (E(1,2) \otimes E(1,2)), \]
\[ \F(E)(3,1)= \Ind_{\s_{2}}^{\s_{3}} (E(2,1) \otimes E(2,1)), \]
\[ \F(E)(2,2) = (E(2,1)\otimes E(1,2))\oplus \Ind_{ \{ 1 \} }
^{\s_{2}\times\s_{2}}(E(1,2) \otimes E(2,1)).\]

Let $(R)$ be the ideal in $\F(E)$ generated
by a right $\s_{3}$-submodule $R(1,3) \subset \F(E)(1,3)$, a left
$\s_{3}$-submodule $R(3,1) \subset \F(E)(3,1)$, and a $(\s_{2},
\s_{2})$-sub-bimodule $R(2,2) \subset \F(E)(2,2)$. Denote by
$\la E;R \ra$ the dioperad $\F(E)/(R)$. A dioperad of
the form $\la E;R \ra$ is called a \emph{quadratic dioperad}
with generators $E$ and relations $R$.
If $\PP$ is a quadratic dioperad, then $\PP(1,1) = \k$
and
$\PP = \la E;R \ra$ where \[ E(i,j) = \PP(i,j) \quad \textrm{for $(i,j) =
(1,2), (2,1)$}, \]
\[ R(i,j) = \Ker(\F(E)(i,j) \rightarrow \PP(i,j)) \quad \textrm{for $(i,j)
= (1,3), (3,1), (2,2)$}. \]

\subsection{}
Let $\sgn_{n}$ denote the sign representation of $\s_{n}$.  
If $V$ is a $(\s_{m}, \s_{n})$-bimodule, then let $V^{*}:=\Hom (V,
\k)$ be the $(\s_{m}, \s_{n})$-bimodule with the transposed
actions, and let $V^{\vee}:= V^{*} \otimes (\sgn_{m} \otimes \sgn_{n})$.
If $\PP = \la E;R \ra$ is a quadratic dioperad, then its
\emph{quadratic dual} is defined by $\PP^{!} := \la E^{\vee};  
R^{\perp}\ra$, where $E^{\vee}(i,j) := E(i,j)^{\vee}$ for   
$(i,j) = (1,2), (2,1)$, and
$R^{\perp}(i,j) \subset \F(E^{\vee})(i,j) = \F(E)(i,j)^{\vee}$ is the
orthogonal complement of $R(i,j)$ for $(i,j) = (1,3), (3,1), (2,2)$.
Here, the identification $\F(E^{\vee})= \F(E)^{\vee}$ is 
canonically defined as follows:
\begin{eqnarray}\label{eqn:2c}
\F(E^{\vee}) = \F(\Lambda^{-1}E^{*}[-1]) 
\stackrel{(\ref{eqn:2a})}{\stackrel{\sim}{\longrightarrow}} 
\Lambda^{-1}\F(E^{*}[-1]) =\F(E)^{\vee}, \end{eqnarray}
where the last equality follows from \cite{GiK} Lemma 3.2.9(b).

{\bf Example.}
The quadratic dioperad $\Bi$ has generators $E(1,2)= \k \cdot l$,
$E(2,1)= \k \cdot \delta$, both with sign actions of $\s_{2}$.
The relations are spanned by
\[ (l_{1}\circ_{1}l)+ (l_{1}\circ_{1}l)\sigma
+(l_{1}\circ_{1}l)\sigma^{2}, \]
\[ (\delta_{1}\circ_{1}\delta)+
\sigma(\delta_{1}\circ_{1}\delta)
+\sigma^{2}(\delta_{1}\circ_{1}\delta), \]
\[ (\delta_{1}\circ_{1}l) - (l_{1}\circ_{1}\delta)
- (l_{1}\circ_{2}\delta) - (l_{2}\circ_{1}\delta) -
(l_{2}\circ_{2}\delta), \]
where $\sigma$ denotes the permutation $(123) \in \s_{3}$.
Thus, $\Bi$-algebras are same as Lie bialgebras \cite{Dr}. 

The quadratic dual $\Bi^{!}$ has generators $E(1,2)= \k \cdot m$,
$E(2,1)= \k \cdot \Delta$, both with trivial actions of $\s_{2}$.
The relations are spanned by
\[ (m_{1}\circ_{1}m)- (m_{1}\circ_{1}m)\sigma, \quad
 (m_{1}\circ_{1}m)- (m_{1}\circ_{1}m)\sigma^{2}, \]
\[ (\Delta_{1}\circ_{1}\Delta)-
\sigma(\Delta_{1}\circ_{1}\Delta), \quad
 (\Delta_{1}\circ_{1}\Delta)-\sigma^{2}(\Delta_{1}\circ_{1}\Delta),\]
\[ (\Delta_{1}\circ_{1}m) - (m_{1}\circ_{1}\Delta), \quad
 (\Delta_{1}\circ_{1}m) - (m_{1}\circ_{2}\Delta), \]
\[ (\Delta_{1}\circ_{1}m) - (m_{2}\circ_{1}\Delta), \quad
 (\Delta_{1}\circ_{1}m) - (m_{2}\circ_{2}\Delta). \]
Thus, a $\Bi^{!}$-algebra is a commutative algebra $A$ equipped with
a cocommutative comultiplication $\Delta : A\rightarrow A\otimes A$ such
that $\Delta$ is a map of $A$-modules. Abrams (\cite{Ab} Theorem 1) has
proved that unital, counital $\Bi^{!}$-algebras are same as unital 
commutative Frobenius algebras.
Observe that $\Bi^{!}(m,n)=\k$ for all $m,n$ (cf. \S\ref{sec:55} below).

{\bf Remarks.}(1) Infinitesimal bialgebras (cf. eg. \cite{Ag}) 
are also algebras over a quadratic dioperad. \\
(2) A dioperad $\PP$ gives a PROP $\QQ$ defined using the same
generators and relations. Algebras over the dioperad $\PP$ are
same as algebras over the PROP $\QQ$. Since the relations defining
$\PP$ are expressed using trees, the grading of a free PROP by 
genus of graphs induces a grading on $\QQ$. The bialgebra PROP
does not have this grading because the compatibility relation
between the multiplication and the comultiplication is expressed
using a tree and a graph of genus $1$ (cf. \cite{ES}), hence
bialgebras cannot be described by a dioperad. Thus, not all
PROP's come from dioperads.

\section{\bf Cobar dual}\label{sec:3}

\subsection{}
Let $\PP$ be a dioperad. The bicomplex $\F(\bar{\PP}^{*}[-1])(m,n)$ 
has differential $d'+d''$, where $d'$ is induced by the differential 
of $\PP$ and $d''$ is induced by edge contractions (\cite{GiK} (3.2.3)):
\begin{eqnarray}\label{eqn:3a}
\bar{\PP}(m,n)^{*} \stackrel{d''}{\rightarrow}
\bigoplus_{|\edge(T)| =1}
\bar{\PP}^{*}(T) \otimes \det(T) \stackrel{d''}{\rightarrow}
\bigoplus_{|\edge(T)| =2}
\bar{\PP}^{*}(T) \otimes \det(T) \stackrel{d''}{\rightarrow} \ldots,
\end{eqnarray}
where the sums are taken over reduced $(m,n)$-trees.
In (\ref{eqn:3a}), the leftmost term is placed in degree $1$.
The \emph{cobar complex} $\C \PP(m,n)$ is defined to be the total complex
of $\F(\bar{\PP}^{*}[-1])(m,n)$. 

Define the \emph{cobar dual} $\D \PP$ of $\PP$ by
$\D\PP:= \Lambda^{-1}\C\PP$. Thus, $\D \PP$ is the total complex of
\begin{eqnarray}\label{eqn:3b}
\bar{\PP}(m,n)^{\vee} \stackrel{d''}{\rightarrow}
\bigoplus_{|\edge(T)| =1}
\bar{\PP}^{*}(T) \otimes \Det(T) \stackrel{d''}{\rightarrow}
\bigoplus_{|\edge(T)| =2}
\bar{\PP}^{*}(T) \otimes \Det(T) \stackrel{d''}{\rightarrow} \ldots,
\end{eqnarray}
where the sums are taken over reduced $(m,n)$-trees. 
In (\ref{eqn:3b}), the leftmost term is placed in degree $3-m-n$.
Note that $\D \PP$ is also the total complex of the free dioperad
generated by $\Lambda^{-1}\bar{P}^{*}[-1]$ with differential $d'+d''$.

\begin{Prop}
There is a canonical quasi-isomorphism $\D\D\PP \rightarrow \PP$.
\end{Prop}
\begin{proof}
(See \cite{GiK} Theorem (3.2.16) or \cite{MSS} Theorem 3.22.)
\[ \D\PP(m,n)  =  \bigoplus_{\textrm{$(m,n)$-trees $T$}}  
\bigotimes_{v\in\Ver(T)} 
\bar{\PP}(\Out(v), \In(v))^{*} \otimes \Det(T).\]
{\setlength\arraycolsep{2pt} \begin{eqnarray*}
\D\D\PP (m,n) & = & \bigoplus_{\textrm{$(m,n)$-trees $T'$}}
\bigotimes_{w\in \Ver(T')}
\big[ \bigoplus_{\textrm{$(\Out(v), \In(v))$-trees $T_{w}$}} 
{} \nonumber\\ & & {}
\bigotimes_{v \in \Ver(T_{w})} 
\bar{\PP}(\Out(v), \In(v))^{*} \otimes \Det(T_{w}) \big]^{*}
\otimes \Det(T')  {} \nonumber\\ & = & {}
\bigoplus_{\textrm{$(m,n)$-trees $T'$}}
\bigoplus_{\textrm{$(m,n)$-trees $T \geq T'$}}
\big(\bigotimes_{v\in\Ver(T)}\bar{\PP}(\Out(v),\In(v))\big)  
{} \nonumber\\ & & {}
\otimes \big(\bigotimes_{w \in \Ver(T')} \Det(T_{w})^{*}\big) 
\otimes \Det(T')
{} \nonumber\\ & = & {}
\bigoplus_{\textrm{$(m,n)$-trees $T$}}
\bar{\PP}(T) \otimes \big(\bigoplus_{T'\leq T} \bigotimes_{w \in \Ver(T')}
\Det(T_{w})^* \otimes \Det(T')\big),
\end{eqnarray*}}
where $T_{w}$ is the subtree of $T$ contracted into $w$.
Note that the summand corresponding to $T$ with just one vertex is
$\PP(m,n)$. 
From here on, continue as in \cite{GiK} p.247.
\end{proof}

\subsection{} 
Let $\PP = \la E;R \ra$ be a quadratic dioperad concentrated in degree
$0$.
Then
\begin{eqnarray}\label{eqn:3k} 
\D\PP(m,n)^{0} = \F(E^{\vee})(m,n), \quad
H^{0}\D\PP(m,n) = \PP^{!}(m,n).
\end{eqnarray}
If the canonical map
$\D\PP(m,n) \rightarrow \PP^{!}(m,n)$ is a quasi-isomorphism for all
$m,n \geq 1$, then $\PP$ is said to be \emph{Koszul}.
If $\PP$ is Koszul, then $\D\PP^{!}$-algebras are called
\emph{strong homotopy $\PP$-algebras}.

\section{\bf Categorical cobar complex}\label{sec:4}
{}From now on, $\PP$ denotes a dioperad concentrated in degree $0$
unless otherwise stated.

\subsection{}
Let $T$ be a forest as in \S\ref{sec:2} (not necessarily reduced). The
direction on edges induces a partial ordering $\preceq$ on $\Ver(T)$. A 
surjective map $\ell:\Ver(T) \rightarrow \{1, \ldots, N\}$ is called a
\emph{$N$-level function} on $T$ if $v \precneqq w$ implies $\ell(v) \lneq
\ell(w)$. A \emph{$N$-level forest} $(T,\ell)$ is a forest $T$ equipped
with a $N$-level function $\ell$. A $N$-level forest  $(T,\ell)$ is said
to be \emph{$i$-saturated} if each path from a leave to a root
traverse some vertex $v$ such that $\ell(v)=i$ ($v$ need not be the same
for all paths). If $(T,\ell)$ is $i$-saturated for all $i$, then we say
that it is \emph{saturated}.

\subsection{}
A \emph{$\s$-bimodule} $\PP=\{\PP(m,n)\}$ is a collection consisting of a
$(\s_{m}, \s_{n})$-bimodule $\PP(m,n)$ for each $m,n \geq 1$. 
We can consider $\k$ as a $\s$-bimodule via
\[ \k(m,n) := \left\{ \begin{array}{ll}
\k & \ \textrm{if $(m,n)=(1,1)$,}\\
0 & \ \textrm{else.} \end{array} \right. \]

Let $\mathcal{M}$ be the category of $\s$-bimodules $\PP$ with 
$\PP(1,1)$ nonzero. Define a monoidal structure on $\mathcal{M}$,
with unit object $\k$, as follows: if $\PP_{1}, \PP_{2}$ are
objects of $\mathcal{M}$, then let
\[ (\PP_{1} \bx \PP_{2})(m,n) :=
\bigoplus_{(T,\ell)}
\bigotimes_{i=1}^{2}
\bigotimes_{v\in\ell^{-1}(i)} \PP_{i}(\Out(v),\In(v)),\]
where the direct sum is taken over all saturated $2$-level 
$(m,n)$-trees $(T,\ell)$.
Thus, if $\PP_{1}, \ldots, \PP_{N}$ are objects of $\mathcal{M}$, then
\begin{eqnarray}\label{eqn:4x}
(\PP_{1} \bx \ldots \bx \PP_{N})(m,n) =
\bigoplus_{(T,\ell)}
\bigotimes_{i=1}^{N}
\bigotimes_{v\in\ell^{-1}(i)} \PP_{i}(\Out(v),\In(v)),
\end{eqnarray}
where the direct sum is taken over all saturated $N$-level 
$(m,n)$-trees $(T,\ell)$. 
Dioperads are precisely monoids in $\mathcal{M}$.

More generally, suppose $\PP_{1}, \ldots, \PP_{N}$ are 
$\s$-bimodules such that $\PP_{i} \in \mathcal{M}$ 
for all $i \leq r$ and $\PP_{i} \notin \mathcal{M}$ 
for all $i > r$. Let $\u{m} = (m_{1}, \ldots, m_{k})$ and
$\u{n} = (n_{1}, \ldots, n_{k})$.
Then we define
$(\PP_{1} \bx \ldots \bx \PP_{N})(\u{m}, \u{n})$ by the right-hand side
of (\ref{eqn:4x}), but with the direct sum taken over all
$N$-level $(\u{m}, \u{n})$-forests $(T,\ell)$ such that $(T,\ell)$ is 
$i$-saturated for all $i \leq r$.

\subsection{} \label{sec:42}
Let $\PP$ be a dioperad. Then the
\emph{categorical cobar complex} $\LL\PP$ of $\PP$ is 
\[ \bar{\PP}^{*} \rightarrow \bar{\PP}^{*\bx 2}
\rightarrow \bar{\PP}^{*\bx 3}
\rightarrow \ldots, \]
where the leftmost term is placed in degree $1$; 
see \cite{MSS} \S3.4 or \cite{SVO} p.398. 
More generally, let $\LL\PP(\u{m},\u{n})$ be the complex
\[ \bar{\PP}^{*}(\u{m},\u{n}) \rightarrow \bar{\PP}^ 
{*\bx 2}(\u{m},\u{n}) \rightarrow \bar{\PP}^{*\bx 3}
(\u{m},\u{n}) \rightarrow \ldots, \]
where the leftmost term is placed in degree $1$.
\begin{Prop} \label{Prop:L}
If $m_{i}+n_{i} \geq 3$ for all $i$, then \[ H(\LL\PP(\u{m},\u{n}))=
H(\C\PP(m_{1},n_{1})) \otimes \ldots \otimes H(\C\PP(m_{k}, n_{k})).\]
In particular, if $m+n\geq 3$, then $H(\LL\PP(m,n))=H(\C\PP(m,n))$.
\end{Prop}

Define a decreasing filtration on $\LL\PP(\u{m},\u{n})$ by letting
$F^{p}\LL\PP(\u{m},\u{n})$ be the subcomplex spanned by 
forests with at least $p$ vertices. 
This gives a spectral sequence
\begin{eqnarray}\label{eqn:4aa}
E_{0}^{p,q} = \frac{F^{p}\LL\PP(\u{m},\u{n})^{p+q}}
{F^{p+1}\LL\PP(\u{m},\u{n})^{p+q}} \Longrightarrow
H^{p+q}(\LL\PP(\u{m},\u{n})). \end{eqnarray}
Hence, Proposition \ref{Prop:L} follows from the following lemma.
\begin{Lem} \label{Lem:L}
If $m_{i}+n_{i} \geq 3$ for all $i$, then
\[ E_{1}^{\bullet,0} \cong \C\PP(m_{1},n_{1}) \otimes \ldots \otimes
\C\PP(m_{k}, n_{k}), \]
and $E_{1}^{p,q} = 0$ if $q \neq 0$.
\end{Lem}

The rest of this section is devoted to the proof of Lemma \ref{Lem:L},
which is same as the proofs in [MSS] \S3.6 and [SVO] \S7 
except for the following modification: instead of using their surjection
algebra, use another algebra that we next define.

\subsection{}\label{sec:43}
Let $X$ be a finite set, and let $S = \{ X_{i} \}_{i \in I}$ be a
collection of nonempty subsets $X_{i}$ of $X$ such that:\\
(a) $X = \bigcup X_{i}$;\\
(b) if $\emptyset \neq Y\subset X_{i}$ for some $i$, then $Y \in S$;\\
(c) if $\emptyset \neq Y\subset X$ and all subsets of $Y$ with $2$ 
elements are in $S$, then $Y \in S$.

Define an associative graded $\k$-algebra $A$ as follows: as a vector
space, $A$ has a basis $\{1\}\cup\{a_{i}\}_{i\in I}$; multiplication of
basis elements is defined by
\[ a_{i}a_{j} := \left\{ \begin{array}{ll}
a_{k} & \ \textrm{if $X_{i}\cap X_{j}=\emptyset$ \& $X_{i} \cup X_{j}=   
  X_{k}$,}\\
0 & \ \textrm{else.} \end{array} \right. \]
The degree of $a_{i}$ is the number of elements of $X_{i}$. Equivalently,
$A$ is a quadratic algebra with generators
$a_{x}$ indexed by elements of $X$, and relations
\[ a_{x}a_{y} := \left\{ \begin{array}{ll}
a_{y}a_{x} & \ \textrm{if $x\neq y$ \& $\{x,y\}\in S$,}\\
0 & \ \textrm{else.} \end{array} \right. \]
Thus, $a_{i}$ is the product of all $a_{x}$ with $x \in X_{i}$.
\begin{Lem}\label{Lem:A}
The quadratic algebra $A$ is Koszul.
\end{Lem}
\begin{proof}
A PBW basis for $A$ is defined by choosing a simple order on the
generators $a_{x}$, and then for each $i$,
choose the admissible expression $a_{i} = a_{x}a_{y}\ldots a_{z}$
to be the one which is maximal in the lexicographical ordering;
cf. \cite{MSS} Proposition 3.65 or \cite{SVO} Proposition 19.
Hence, $A$ is Koszul by \cite{Pr} Theorem 5.3.
\end{proof}

\subsection{} \label{sec:44}
Fix a forest $T$ such that each component of $T$ has at least one vertex. 
Call $V \subset \Ver(T)$ a \emph{level} of $T$ if there
exists a level function $\ell$ with $V = \ell^{-1}(i)$ for some $i$.   
Let $X=\Ver(T)$ and let $S$ be the collection of all levels of $T$. 
Then $X$ and $S$ satisfy the conditions in \S \ref{sec:43}.
(To verify condition (c), use induction on the number of elements in $Y$.) 
Thus, there is a corresponding algebra $A$. Recall that $A$ has a basis
labelled by $1$ and all the levels.

Call a sequence $V_{1}, \ldots, V_{N}$ of levels \emph{good} if there
exists a $N$-level function $\ell$ such that $V_{i} = \ell^{-1}(i)$ for
all $i$. Consider the cochain complex $C^{\bullet}(A) :=
\bar{A}^{*\otimes \bullet}$. For each $N$,
there is a direct sum decompositon $\bar{A}^{\otimes N} = U_{N} \oplus
V_{N}$, where $U_{N}$ is spanned by $a_{i_{1}}\otimes \ldots \otimes
a_{i_{N}}$ such that $X_{i_{1}}, \ldots, X_{i_{N}}$ is not a good sequence
of levels, and $V_{N}$ is spanned by $a_{i_{1}}\otimes \ldots \otimes
a_{i_{N}}$ such that $X_{i_{1}}, \ldots, X_{i_{N}}$ 
is a good sequence of levels. 
Let $C_{U}^{N}$ be the annihilator of $V_{N}$ and $C_{V}^{N}$ be the
annihilator of $U_{N}$.
Denote by $\partial$ the boundary map of the chain complex
$\bar{A}^{\otimes \bullet}$. Since $\partial(U_{N}) \subset U_{N-1}$ and
$\partial(V_{N}) \subset V_{N-1}$, it follows that
$C^{\bullet}(A) = C_{U}^{\bullet} \oplus   
C_{V}^{\bullet}$, a direct sum of two subcomplexes.

Let $\{a_{i}^{*}\}$ be the dual basis of $\{a_{i}\}$, then $C_{V}^{N}$ has
a basis consisting of all $a_{i_{1}}^{*}\otimes \ldots \otimes
a_{i_{N}}^{*}$ for which
$X_{i_{1}}, \ldots, X_{i_{N}}$ is a good sequence of levels. The
differential on $C_{V}^{N}$ is
\[ \delta(a_{i_{1}}^{*}\otimes \ldots \otimes a_{i_{N}}^{*}) = \sum
(-1)^{r} a_{i_{1}}^{*}\otimes \ldots \otimes \delta(a_{i_{r}}^{*})     
\otimes \ldots \otimes a_{i_{N}}^{*}, \]
\[ \delta(a_{i}^{*}) = \sum_{X_{j} \coprod X_{k} = X_{i}}
a_{j}^{*} \otimes a_{k}^{*}  \]

\begin{Lem}\label{Lem:T}
If $N \neq |\Ver(T)|$, then $H^{N}(C_{V}^{\bullet})= 0$.
If $N =|\Ver(T)|$, then $H^{N}(C_{V}^{\bullet}) = \det(T)$.
\end{Lem}
\begin{proof}
(See \cite{MSS} Proposition 3.67 or \cite{SVO} Proposition 21.)

Note that $H(C^{\bullet}(A)) = H(C_{U}^{\bullet}) \oplus
H(C_{V}^{\bullet})$.
Let $[a_{i_{1}}^{*}\otimes \ldots \otimes a_{i_{N}}^{*}] \in  
H^{N}(C_{V}^{\bullet})$ be a nonzero cohomology class. 
By Lemma \ref{Lem:A}, $A$ is Koszul. Therefore,
$X_{i_{1}}, \ldots , X_{i_{N}}$ must each be a one element  
subset of $\Ver(T)$. But $X_{i_{1}}, \ldots , X_{i_{N}}$ is a good
sequence of levels, hence $N=|\Ver(T)|$.

Define a map $\iota: C_{V}^{N} \rightarrow \bigwedge^{|\Ver(T)|}
\k^{\Ver(T)} \cong \det(T)$ by sending $a_{x}^{*}\otimes a_{y}^{*}\otimes
\ldots \otimes a_{z}^{*}$ to $x \wedge y \wedge \ldots \wedge z$, where
$\{x\}, \{y\},\ldots, \{z\}$ is any good sequence of levels.
The coboundaries in $C_{V}^{N}$ are  
spanned by terms of the form
$a_{x}^{*}\otimes \ldots \otimes \delta a_{ \{ y_{1}, y_{2} \} }^{*}
\otimes \ldots \otimes a_{z}^{*}$ with 
\[ \delta a_{ \{ y_{1}, y_{2} \} }^{*} =
a_{y_{1}}^{*} \otimes a_{y_{2}}^{*}+ a_{y_{2}}^{*} 
\otimes a_{y_{1}}^{*},\]
hence $\iota$ induces a map $H^{N}(C_{V}^{\bullet}) \rightarrow \det(T)$.
Clearly, $\iota$ is surjective. Observe that any two good sequences of
singleton levels can be obtained from each other by performing
transpositions of succesive levels while remaining good at the
intermediate steps.
This implies that $H^{N}(C_{V}^{\bullet})$ is at most one dimensional.
Hence $\iota$ is bijective.
\end{proof}
 
\subsection{}
Suppose $m_{i}+n_{i} \geq 3$ for all $i$.
For a $(\u{m},\u{n})$-forest $T$ with $p$ vertices, let
\[ E_{0}^{p,\bullet}(T) := \bar{\PP}^{*}(T_{1})\otimes \cdots 
\otimes \bar{\PP}^{*}(T_{k}) \otimes C_{V}^{p+\bullet},\]
where $T_{1}, \ldots, T_{k}$ are the connected components of $T$.
Observe that the complex $E_{0}^{p, \bullet}$ of (\ref{eqn:4aa}) is the
direct sum of $E_{0}^{p,\bullet}(T)$ over all $(\u{m},\u{n})$-forests 
$T$ with $p$ vertices. Therefore, Lemma \ref{Lem:T} implies 
Lemma \ref{Lem:L}.

\section{\bf Koszul complex}\label{sec:5}
{}From now on, assume that $\PP$ is quadratic unless otherwise stated.

\subsection{}
Observe that a summand in $\PP \bx (\Lambda \PP^{!})^{*}$
corresponding to a saturated $2$-level tree $(T, \ell)$ is placed in 
degree $\sum_{v\in \ell^{-1}(2)}(2- |\Out(v)| - |\In(v)|)$.
The \emph{Koszul complex} of a quadratic dioperad $\PP = \la E;R \ra$ is
defined to be $\K\PP := \PP \bx (\Lambda \PP^{!})^{*}$ with
differential $d$ defined as follows.

Let \[ \mu:\PP\bx\PP \rightarrow \PP \quad \textrm{and}\quad
\delta: (\Lambda \PP^{!})^{*} \rightarrow
(\Lambda \PP^{!})^{*} \bx (\Lambda \PP^{!})^{*} \]
be, respectively, the composition map of $\PP$ and the 
dual of the composition map of $\Lambda\PP^{!}$.
Let \[ \alpha(1,2): E(1,2)[1] \rightarrow E(1,2), \qquad
\alpha(2,1): E(2,1)[1] \rightarrow E(2,1) \]
be the identity maps of degree $1$. 
Then $d$ is defined to be the composition
\[ \PP \bx (\Lambda \PP^{!})^{*}
\stackrel{\Id \bx \delta}{\longrightarrow}
\PP \bx (\Lambda \PP^{!})^{*} \bx (\Lambda \PP^{!})^{*}
\stackrel{\Phi}{\longrightarrow}
\PP \bx \PP \bx (\Lambda \PP^{!})^{*}
\stackrel{\mu\bx\Id}{\longrightarrow}
\PP \bx (\Lambda \PP^{!})^{*}, \]
where $\Phi$ is defined to be $\Id\bx\alpha\bx\Id$ on saturated
$3$-level trees whose second level has one vertex of
valency $3$ and all other second level vertices are of valency $2$,
and zero otherwise. 

The proof of the next lemma is a standard argument for twisted
tensor products.
\begin{Lem} The map $d$ is a differential, i.e. $dd=0$. \end{Lem}
\begin{proof}
{\setlength\arraycolsep{2pt}
\begin{eqnarray*}
dd & = & (\mu\bx\Id) \Phi (\Id\bx\delta) (\mu\bx\Id) \Phi (\Id\bx\delta)
\nonumber\\ & = &
(\mu\bx\Id) \Phi (\mu\bx\Id\bx\Id) (\Id\bx\Id\bx\delta)  \Phi
(\Id\bx\delta) \nonumber\\ & = &
(\mu\bx\Id) (\mu\bx\Id\bx\Id) (\Id\bx\Phi) (\Phi\bx\Id)
(\Id\bx\Id\bx\delta) (\Id\bx\delta) \nonumber\\ & = &  
(\mu\bx\Id) (\Id\bx\mu\bx\Id) (\Id\bx\Phi) (\Phi\bx\Id)
(\Id\bx\delta\bx\Id) (\Id\bx\delta). 
\end{eqnarray*}}
But $(\Id\bx\mu\bx\Id) (\Id\bx\Phi) (\Phi\bx\Id)
(\Id\bx\delta\bx\Id) = 0$ by relations defining $\PP$ (in the case when 
the second level has one vertex of valency $4$ and all other 
second level vertices are of valency $2$)
and by Koszul sign rule (in the case when the second level has
two vertices of valency 3 and all other second level vertices are of
valency $2$).
\end{proof}

\begin{Thm}\label{Thm:K}
A quadratic dioperad $\PP$ is Koszul if and only if
the canonical map $\K\PP\rightarrow\k$ is a quasi-isomorphism
(i.e. $\K\PP(m,n)$ is exact for all $(m,n) \neq (1,1)$).
\end{Thm}
\begin{proof}
As usual, there is an exact sequence $\B\PP(m,n)$:
\[ (\PP\bx\bar{\PP}\bx\cdots\bx\bar{\PP})(m,n) \rightarrow \cdots
\rightarrow (\PP\bx\bar{\PP})(m,n) \rightarrow 
\PP(m,n) \rightarrow \k(m,n), \]
where the leftmost term is placed in degree $1-m-n$ and the rightmost term 
is placed in degree $0$. Suppose $m+n \geq 3$.
Define a decreasing filtration on $\B\PP(m,n)$ by letting 
$F^{p}\B\PP(m,n)$ be the subcomplex spanned by  
$1$-saturated $(m,n)$-trees $(T, \ell)$ with
\[\sum_{v \in \ell^{-1}(1)}(|\Out(v)|+|\In(v)|-2) \geq p.\]
This gives a spectral sequence
\begin{eqnarray}\label{eqn:5a}
E_{0}^{p,q} = \frac{F^{p}\B\PP(m,n)^{p+q}}
{F^{p+1}\B\PP(m,n)^{p+q}} \Longrightarrow H^{p+q}(\B\PP(m,n))=0 .
\end{eqnarray}

Let $(T,\ell)$ be a saturated $2$-level $(m,n)$-tree. 
Suppose
\[ \{v_{1},\ldots, v_{k}\} := \{v \in \ell^{-1}(2);\ |\Out(v)|+ 
|\In(v)| \geq3 \} \neq \emptyset.\] Let
\[ \u{r}:=(|\Out(v_{1})|, \ldots, |\Out(v_{k})|) \quad\textrm{and}
\quad \u{s}:=(|\In(v_{1})|, \ldots, |\In(v_{k})|).\]
Note that, by (\ref{eqn:2f}), we have 
{\setlength\arraycolsep{2pt}
\begin{eqnarray}\label{eqn:5b}
|\u{r}|+|\u{s}|-2k + 
\sum_{v \in \ell^{-1}(1)}(|\Out(v)|+|\In(v)|-2)
& = & m+n-2. 
\end{eqnarray}}
Let 
\begin{eqnarray}\label{eqn:5t}
E_{0}^{p,\bullet}(T,\ell) := 
\bigotimes_{v\in\ell^{-1}(1)} \PP(\Out(v),\In(v))
\otimes \LL\PP(\u{r},\u{s})^{*}[p+1]. \end{eqnarray}
Observe that, for $0 \leq p \leq m+n-3$,
the complex $E_{0}^{p, \bullet}$ of (\ref{eqn:5a}) is the
direct sum of $E_{0}^{p,\bullet}(T,\ell)$ over all saturated 
$2$-level $(m,n)$-trees $(T,\ell)$ such that
\begin{eqnarray}\label{eqn:5c} 
\sum_{v \in \ell^{-1}(1)}(|\Out(v)|+|\In(v)|-2) =p .
\end{eqnarray}
By (\ref{eqn:3k}) and Proposition \ref{Prop:L},
$H^{|\u{r}|+|\u{s}|-2k} (\LL\PP(\u{r},\u{s}))
= (\Lambda\PP^{!})(\u{r},\u{s})$. Hence, by 
(\ref{eqn:5b})--(\ref{eqn:5c}),
\begin{eqnarray}\label{eqn:5k}
E_{1}^{\bullet, 1-m-n} = \K\PP(m,n)[2-m-n].
\end{eqnarray}

Now assume that $\PP$ is Koszul. Suppose $(m,n) \neq (1,1)$. 
By Proposition \ref{Prop:L} and (\ref{eqn:5b})--(\ref{eqn:5c}),
$E_{1}^{p,q}=0$ if $q\neq 1-m-n$. Hence, by (\ref{eqn:5k}),
$\K\PP(m,n)$ is exact.

Conversely, assume $\K\PP(m,n)$ is exact for all $(m,n) \neq (1,1)$.
Let $\u{w}:=(w_{1},\ldots,w_{j})$ and $\u{x}:=(x_{1},\ldots,x_{j})$,
where $w_{l}+ x_{l} \geq 3$ for all $l$.
We shall prove, by induction on $|\u{w}|+|\u{x}|-2j$, that 
\begin{eqnarray}\label{eqn:5x}
H^{i}(\LL\PP(\u{w},\u{x}))=0 \quad \textrm{if}\ i< |\u{w}|+|\u{x}|-2j.
\end{eqnarray}
By Proposition \ref{Prop:L}, this would imply that $\PP$ is Koszul.

First, note that (\ref{eqn:5x}) is true if
$|\u{w}|+|\u{x}|-2j = 1$.
Assume that (\ref{eqn:5x}) is true if $|\u{w}|+|\u{x}|-2j<N$.
Now let $|\u{w}|+|\u{x}|-2j=N$. By Proposition \ref{Prop:L}, 
it suffices to show that, for each $l$, we have 
$H^{i}(\LL\PP(w_{l},x_{l}))=0$ if $i< w_{l}+x_{l}-2$.
Consider the spectral sequence (\ref{eqn:5a}) with 
$m:=w_{l}$ and $n:= x_{l}$. 
Note that by (\ref{eqn:5t})--(\ref{eqn:5c}),
$E_{0}^{0, \bullet}= \LL\PP(m,n)^{*}[1]$,
so $E_{1}^{0,q}= H^{q+1}(\LL\PP(m,n)^{*})$.
Moreover, by induction hypothesis and (\ref{eqn:5b})--(\ref{eqn:5c}), 
we have $E_{1}^{p,q}=0$ if $p>0$ and $q > 1-m-n$.
Hence, by (\ref{eqn:5k}) and exactness of $\K\PP(m,n)$,
we obtain
\[ E_{2}^{p,q} = \left\{ \begin{array}{ll}
H^{q+1}(\LL\PP(m,n)^{*}) & \ \textrm{if $p=0$ and $q\geq 2-m-n$,}\\
0 & \ \textrm{else.} \end{array} \right. \]
Therefore, $H^{i}(\LL\PP(m,n))=0$ for $i<m+n-2$.
\end{proof}

\subsection{}\label{sec:55}
Observe that if $\PP$ is a quadratic dioperad, then 
$\A(n):= \PP(1,n)$ and $\BB(n):= \PP^{\op}(1,n)$ are 
quadratic operads. The inclusions $\A \hookrightarrow \PP$
and $\BB^{\op} \hookrightarrow \PP$ induce canonical morphisms
$\A \bx \BB^{\op} \rightarrow \PP$ and $\BB^{\op} \bx \A \rightarrow \PP$.
The following proposition is adapted from \cite{M2}; cf. also \cite{FM}.
\begin{Prop}\label{Prop:M}
Let $\PP= \la E;R \ra$ be a quadratic dioperad,
let $\A(n):= \PP(1,n)$, and let $\BB(n):= \PP^{\op}(1,n)$.\\
\emph{(a)} If $\PP(m,n) = (\A \bx \BB^{\op})(m,n)$ (resp. $\PP(m,n) =
(\BB^{\op} \bx \A)(m,n)$) for $(m,n) = (2,2)$, $(2,3)$, and $(3,2)$, then
$\PP= \A \bx \BB^{\op}$ (resp. $\PP= \BB^{\op} \bx \A$).\\
\emph{(b)} We have $\PP= \A \bx \BB^{\op}$ if and only if
$\PP^{!}= \BB^{! \op} \bx \A^{!}$.\\
\emph{(c)} If $\A$ and $\BB$ are Koszul, and $\PP= \A \bx \BB^{\op}$,
then $\PP$ is Koszul. 
\end{Prop}
\begin{proof}
(a) (See \cite{M2} Theorem 2.3.)

Note that $\F(E)(2,2)= (\A \bx \BB^{\op})(2,2) \oplus
(\BB^{\op} \bx \A)(2,2)$. We have $\PP(2,2)= (\A \bx \BB^{\op})(2,2)$ 
if and only if there is a $(\s_{2}, \s_{2})$-equivariant map
\[ \lambda: (\BB^{\op} \bx \A)(2,2) \rightarrow (\A \bx \BB^{\op})(2,2)
\] such that
\begin{eqnarray*}
R(2,2)= \{ x-\lambda(x);\ x\in (\BB^{\op} \bx \A)(2,2) \}.
\end{eqnarray*}
The rest of the proof is same as \cite{M2} Theorem 2.3.
Similarly in the other case.\\
(b) (See \cite{M2} Lemma 4.3.)

Assume that $\PP= \A \bx \BB^{\op}$, and let $\lambda$ be 
the map in the proof of (a). Via the identification (\ref{eqn:2c}),
we obtain
\begin{eqnarray*}
R^{\perp}(2,2)= \{ \alpha- \lambda^{\vee}(\alpha);\ \alpha\in
(\A^{!} \bx \BB^{!\op})(2,2)\}, \end{eqnarray*}
where $\lambda^{\vee}$ is the dual of $\lambda$.
It follows by a direct verification, using (a), that
$\PP^{!}= \BB^{!\op} \bx \A^{!}$. Similarly for the converse.\\
(c) (See \cite{M2} Proposition 4.4 and Theorem 4.5.)

Observe that $\bar{\PP} \bx \bar{\PP} \subset \F(\bar{\PP})$.
Thus, using (\ref{eqn:2a}), we have
\[ \K\PP = \A \bx \BB^{\op} \bx (\Lambda(\BB^{!\op} \bx \A^{!}))^{*}
= \A \bx \BB^{\op} \bx (\Lambda\BB^{\op!})^{*} \bx (\Lambda\A^{!})^{*}.
\]
Suppose $(m,n)\neq (1,1)$. Define a decreasing filtration on
$\K\PP(m,n)$ by letting  $F^{p}\K\PP(m,n)$ be the subcomplex 
spanned by saturated $4$-level $(m,n)$-trees $(T,\ell)$ with
\[ \sum_{v\in \ell^{-1}(4)} (2- |\Out(v)|- |\In(v)|) \geq p.\]
This gives a spectral sequence
\[ E_{0}^{p,q} = \frac{F^{p}\K\PP(m,n)^{p+q}}
{F^{p+1}\K\PP(m,n)^{p+q}} \Longrightarrow
H^{p+q}\K\PP(m,n). \]
Now $E_{0}^{p, \bullet}$ is the subcomplex of
$\A \bx \K\BB^{\op} \bx (\Lambda\A^{!})^{*}$ (with differential
$\Id \bx d \bx \Id$) spanned by 
saturated $3$-level $(m,n)$-trees $(T,\ell)$ with 
\[ \sum_{v\in \ell^{-1}(3)} (2- |\Out(v)|- |\In(v)|) =p.\]
Thus, $E_{1}^{p,q}=0$ if $q \neq 0$ and 
$E_{1}^{\bullet, 0} = \K\A(m,n)$. Hence, $E_{2}^{p,q} = 0$ for all
$p,q$.
\end{proof}

\begin{Cor}
The dioperad $\Bi$ is Koszul.
\end{Cor}
\begin{proof}
This is immediate from Proposition \ref{Prop:M} and Koszulity of the 
Lie operad $\mathscr{L}ie$ (\cite{BG} (6.4) or \cite{GiK} Corollary
(4.2.7)). We have: 
$\Bi = \mathscr{L}ie \bx \mathscr{L}ie^{\op}$,
$\Bi^{!}= \mathscr{C}om^{\op} \bx \mathscr{C}om$, where
$\mathscr{C}om = \mathscr{L}ie^{!}$.
\end{proof}

{\bf Remarks.}
(1) Note that, in more explicit terms, the cobar dual of $\Bi^{!}$ is the
complex of trees with multiple leaves and multiple roots. 
The Koszulity of $\Bi^{!}$ means that the cohomology of
its cobar dual is zero in all negative degrees, but the author does not
know of a direct proof of this.\\
(2) Similarly, one can show that the infinitesimal
bialgebra dioperad is Koszul.

\footnotesize{

\section*{\bf Acknowledgements}
I am very grateful to Victor Ginzburg for suggesting the problem and
for many useful discussions. I also thank Martin Markl
for answering questions on his paper \cite{M2}.

}

\end{document}